\documentstyle[12pt]{article}
\def\baselinestretch{1.0}
\title
{\bf  Hausdorff dimension,\\ Mean quadratic variation of \\infinite self-similar
measures\footnotemark }
\author
{  Zu-Guo Yu$^1$, Fu-Yao Ren$^2$ and Jin-Rong Liang$^2$\\
 {\small $^1$Institute of Theoretical Physics, Academia Sinica},\\
  {\small P.O. Box 2735, Beijing 100080, P.R.C. E-mail:yuzg@itp.ac.cn.}\\
 {\small $^2$Institute of Mathematics, Fudan University , Shanghai 200433, 
 P.R.C.}\\
\\
 {\small (Received on 10 October, 1996; Accepted on 15 August, 1997)}
 }
\date{}
\newcommand{\be}{\begin{equation}}
\newcommand{\ee}{\end{equation}}
\begin{document}
\newtheorem{Theorem}{\quad Theorem}
\newtheorem{Proposition}{\quad Proposition}
\newtheorem{Definition}{\quad Definition}
\newtheorem{Lemma}{\quad Lemma}
\newtheorem{Corollary}{\quad Corollary}
\maketitle
 \renewcommand{\thefootnote}{\fnsymbol{footnote}}
 \footnotetext{* Project partially supported by the Tianyuan Foundation of China and
Postdoctoral station Founda-tion of the State Education Committee. } 
\vskip 10mm

{\bf Abstract:} \ \ Under weaker condition than that of Riedi \& Mandelbrot, the 
Hausdorff (and
Hausdorff-Besicovitch) dimension of infinite self-similar set $K\subset {\bf R}^d$ which 
is the invariant compact set of
 infinite contractive similarities $\{S_j(x)=\rho_jR_jx+b_j\}_{j\in {\bf N}}$
($0<\rho_j<1,b_j\in {\bf R}^d, R_j$ orthogonal) satisfying open set
condition is obtained.
It is proved (under some additional hypotheses) that the $\beta$-mean quadratic variation 
of
  infinite self-similar measure is of
asymptotic property (as $t\longrightarrow 0$).
\vskip 5mm

\vskip 5mm
\pagebreak

\section{Introduction}
\ \ \ \   In this paper, we denote ${\bf R}^d$ the $d$-dimensional Euclidean space,
${\bf N}$ the set of natural numbers and ${\bf Z}$ the set of integer numbers.

  For given finite contractive similarities $\{S_j(x)=\rho_jR_jx+b_j\}_{j=1}^{m}$ of ${\bf R}^d$,
where $0<\rho_j<1,b_j\in {\bf R}^d, R_j$ orthogonal, J.E.Hutchinson [1]
proved that there exists unique compact set
$K_1$ satisfying
$$K_1=\cup_{j=1}^{m}S_j(K_1).$$
 $K_1$ is called
{\it self-similar set}. If there exists an open set $O_1$ satisfying
$S_j(O_1)\subset O_1$ and $S_i(O_1)\cap S_j(O_1)=\emptyset$ $(i\neq j)$, we call
that $\{S_j\}_{j=1}^{m}$ satisfy {\it open set condition}. We call that they
satisfy {\it strong open set condition} if the sets $S_j(\overline{O})$ are disjoint.
 Then

{\bf Theorem A} (Hutchinson)
\ \ \ \ 
 If $\{S_j\}_{j=1}^m$ satisfy open set condition, then the Hausdorff dimension $s'$
  of $K_1$
 is the unique solution of the equation  $\sum_{j=1}^{m}
\rho_j^{s'}=1$.

 In [1], he also proved that for given probability vector
$P=(P_1,P_2,\cdots,P_m)$ satisfying $\sum_{j=1}^{m}P_j=1$, there exists unique
probability measure $\mu_1$ on ${\bf R}^d$ satisfying
$$\mu_1(\cdot)=\sum_{j=1}^{m}P_j\mu_1(S_j(\cdot))$$
and the support set of $\mu_1$ is $K_1$. $\mu_1$ is called {\it self-similar
measure} and $\{P_j\}_{j=1}^m$ is called {\it weights} of $\mu_1$.

  Ka-Sing Lau and Jian-rong Wang [2], and R.S.Strichartz [3-7] have
done much study on Fourier analysis of self-similar measure. R.S.Strichartz in
[3] (or [7]) discussed many fractal measures. 
  If $\mu$ is self-similar measure on ${\bf R}^d$, Strichartz [4-6] discussed
the asymptotic property (as $r\longrightarrow\infty$) of function
$$H(r)=\frac{1}{r^{d-\beta'}}\int_{|x|\leq r}|F(x)|^2dx,$$
where $F(x)=(d\mu)^{\wedge}$ and $\beta'$ is defined by $\sum_{j=1}^{m}\rho_j^{-\beta'}P_j^2=1$.

  Let $\mu$ be a Borel measure on ${\bf R}^d$, $f$ be a Borel measurable function,
we use $\mu_f$ to denote the measure defined by $\mu_f(E)=\int_Efd\mu$ for any
Borel set $E$ in ${\bf R}^d$.

  It is proved in [5] that if $\{S_j\}_{j=1}^m$ satisfies the strong open set condition,
then for the self-similar measure $\mu$ defined by natural weights (i.e. $P_j=\rho_j^{\beta'}, \quad \beta'=s'$)
$$\frac{1}{r^{d-\beta'}}\int_{|x|\leq r}|(\mu_f)^{\wedge}|^2dx=q(r)\int|f|^2d\mu+
E(r)\quad \hbox{ for}\quad \forall f\in L^2(d\mu), \qquad \qquad (*)$$
where $E(r)\longrightarrow 0$ as $r\longrightarrow +\infty$, and $q(r)$ is a multiplicative
periodic function or a positive constant.

  Let $\mu$ be a $\sigma$-finite measure on ${\bf R}^d$, for $0\leq \alpha \leq d$,
let $$V_{\alpha}(t;\mu)=\frac{1}{t^{d+\alpha}}\int_{{\bf R}^d}|\mu(B_t(x))|^2dx,$$
where $B_t(x)$ is the ball of radius $t$, centered at $x$. We will call
$\lim\sup_{t\longrightarrow 0}V_{\alpha}(t;\mu)$ the {\it upper $\alpha$-mean
quadratic variation} (m.q.v.) of $\mu$, and simply call it $\alpha$-m.q.v. if the
limit exists.

  If $\mu$ is a self-similar measure on ${\bf R}^d$, Ka-sing Lau and Jian-rong Wang
[2] proved the following two Theorems

{\bf Theorem B}([2])\ \ \  Under some additional conditions, we have
$$\lim_{ t\longrightarrow 0}[V_{\beta'}(t;\mu)-p(t)]=0.$$
where $p(t)$ is a multiplicative
periodic function or a positive constant and $\beta'$ is defined as above.

{\bf Theorem C} ([2])
\ \ \ 
 If the self-similar
 measure $\mu$ defined by natural weights (i.e. $P_j=\rho_j^{\beta'}, \beta'=s'$), under some additional
 hypotheses
$$\lim_{t\longrightarrow 0}[\frac{1}{t^{d+\beta'}}\int_{{\bf R}^d}|\mu_f(B_t(x))|^2dx-p(t)
\int|f|^2d\mu]=0\quad \hbox{ for}\quad \forall f\in L^2(d\mu), $$
where $p(t)$ is the function in Theorem B.

  R.H.Riedi and B.B.Mandelbrot [8] introduced infinite self-similar sets and
infinite self-similar measures on ${\bf R}^d$ (definitions see later of this paper),
 discussed multifractal formalism for infinite self-similar measures and
 the Hausdorff dimension of infinite
self-similar sets (under some additional conditions). In this paper, under weaker
condition than that of Riedi \& Mandelbrot, we extend Theorem A to the infinite
self-similar case.
 If $\mu$ is infinite self-similar measure
and the equation $\sum_{j=1}^{\infty}P_j^2\rho_j^{-\beta}=1$ has finite solution $\beta$,
then under some additional hypotheses, R.S.Strichartz [5] obtained the asymptotic
property of function $H(r)$ and conclusion (*).
In this paper, we also extend Theorem B,C to the infinite self-similar case.

\section { Hausdorff (and Hausdorff-Besicovitch) dimension of infinite self-similar set.}
 \ \ \ \   For given infinite contractive similarities $\{S_j(x)=\rho_jR_jx+b_j\}_{j\in {\bf N}}$ of ${\bf R}^d$,
where $0<\rho_j<1,b_j\in {\bf R}^d, R_j$ orthogonal,
from [8], there exists unique compact set $K$ satisfying
$$K=\overline{\cup_{j=1}^{\infty}S_j(K)}.$$
$K$ is called {\it infinite self-similar set}. $K$ can be constructed as following.
Let $E_0\subset {\bf R}^d$ be a compact set, denote $E_{j_1\cdots j_k}=S_{j_1}\circ\cdots\circ S_{j_k}(E_0)$,
then $$K=\cap_{k=0}^{\infty}\overline{\cup_{j_1,\cdots,j_k\in {\bf N}}E_{j_1\cdots j_k}}.$$
For given probability sequence $(P_1,P_2,\cdots)$ with $\sum_{j=1}^{\infty}P_j=1$, from
[8], there exists unique probability measure $\mu$ on ${\bf R}^d$ satisfying
$$\mu(\cdot)=\sum_{j=1}^{\infty}P_j\mu(S_j(\cdot)).$$
We call $\mu$ {\it infinite self-similar measure} and $\{P_j\}_{j=1}^{\infty}$
{\it weights } of $\mu$.
Its support set is $K$.

\begin{Definition} We call $\{S_j(x)\}_{j\in {\bf N}}$ satisfying {\it open set
condition} if there exists a bounded open set $O\subset$ ${\bf R}^d$ such that
$S_j(O)\subset O$ and $S_i(O)\cap S_j(O)=\emptyset$ $(i\neq j)$.
\end{Definition}

  For any subset $A \in {\bf R}^d$  and $0\leq s<\infty$, let ${\cal M}_{\delta}^s(A)=
  \inf\sum_{i=1}^{\infty}|A_i|^s$,
where $A=\cup_{i=1}^{\infty}A_i$ is a countable decomposition of $A$ into subsets
of diameter $|A_i|<\delta$ ($>0$). We set $|A_i|^0=0$ if $A_i$ is empty and $|A_i|^0=1$ otherwise.
The the $s$-dimensional measure of $A$ is defined to be
$${\cal M}^s(A)=\sup_{\delta>0}{\cal M}_{\delta}^s(A).$$
The Hausdorff-Besicovitch dimension$^{[9]}$ of $A$ is
$$\dim_M(A)=\sup\{0\leq s<\infty:\quad {\cal M}^s(A)>0\}.$$

  {\bf Remark}: It is easy to see that in the definition of ${\cal M}_{\delta}^s(A)$,
we can replace $|A_i|$ by $|\overline{A_i}|$.

  From the definition of fractal dimension$^{[10]}$ $\dim_H(A)$, we can see that
\be \dim_H(A)\leq \dim_M(A).  \ee

\begin{Theorem} If the equation $\sum_{j=1}^{\infty}\rho_j^s=1$ has finite
solution $s$, and $\{S_j\}_{j=1}^{\infty}$ satisfy open set condition, $K$
 is the infinite self-similar set,
  then the
Hausdorff-Besicovitch dimension $\dim_M(K)$ and Hausdorff dimension $\dim_H(K)$
of $K$ is $s$.
\end{Theorem}

  {\bf Remark}. Our condition is weaker than Riedi \& Mandelbrot's [8] condition:
there exist numbers $r$, $R$ such that $-\infty<\log r\leq (1/j)\log \rho_j\leq \log R<0$ $\forall j$.

{\bf Proof of Theorem 1} To get the upper bound. Let $K=\cup_{i=1}^{\infty}A_i$ be any
decomposition of $K$ into subsets of diameter $<\delta$, then a new
decomposition is provided by $K=\cup_{i=1}^{\infty}\overline{\cup_{j=1}^{\infty}A_{ij}}$,
where $A_{ij}=\varphi_j(A_i)$. Because
\begin{eqnarray*} 
\sum_{i=1}^{\infty}\sum_{j=1}^{\infty}|A_{ij}|^s&\leq& \sum_{i=1}^{\infty}\sum_{j=1}^{\infty}|\rho_j|^s|A_i|^s\\
                   &\leq& (\sum_{j=1}^{\infty}\rho_j^s)\sum_{i=1}^{\infty}|A_i|^s,
                   \end{eqnarray*}
it follows that whenever $\sum_{j=1}^{\infty}\rho_j^s<1$ we must have ${\cal M}_{\delta}^s(K)=0$,
then ${\cal M}^s(K)=0$. As $\dim_M(K)=\inf\{s:\quad {\cal M}^s(K)=0\}$, hence
$\dim_M(K)\leq s$ where $\sum_{j=1}^{\infty}\rho_j^s=1$. From (1), we have $\dim_H(K)\leq s$.

  To get the lower bound. We let $K^{(m)}$ be the self-similar set generated by
$\{S_j\}_{j=1}^m$, then from Theorem 8 of ref.[11], we have
\be \dim_M(K^{(m)})\geq \min\{d,s^{(m)}\}, \ee
where $s^{(m)}$ is the positive solution of $\sum_{j=1}^m\rho_j^{s^{(m)}}=1$. Using
Theorem 4.13 of ref.[10], similar to the proof of Theorem 8 of ref.[11], we can
obtain
\be \dim_H(K^{(m)})\geq \min\{d,s^{(m)}\}. \ee
Then from Lemma 8 of ref.[8], we have $\lim_{m\longrightarrow\infty}s^{(m)}=s$,
where $\sum_{j=1}^{\infty}\rho_j^s=1$. Since for any $m$, $K^{(m)}\subset K$,
we have $\dim_M(K)\geq\dim_M(K^{(m)})$ and $\dim_H(K)\geq\dim_H(K^{(m)})$.
From open set condition, we have $s<d$, then from (2) and (3), we have
\be \dim_M(K)\geq s^{(m)} \ee
and
\be \dim_H(K)\geq s^{(m)}. \ee
Take limit from (4) and (5), we have $\dim_M(K)\geq s$ and $\dim_H(K)\geq s$.
\ \ \ \#

  The method used in proof of Theorem 1 can be used to estimate the Hausdorff
(and Hausdorff-Besicovitch) dimension of the limit set of infinite non-similar
contractive maps.

\begin{Corollary} Let $\{\varphi_j\}_{j=1}^{\infty}$ be infinite contractive
maps with
$$|\varphi_j(x)-\varphi_j(y)|\leq c_j |x-y|,\quad x,y\in  {\bf R}^d,\quad j=1,2,\cdots,$$
and satisfying open set condition, and denote $E$ their contractive-invariant set.
  If the equation $\sum_{j=1}^{\infty}c_j^u=1$
has finite solution $u$, then $\dim_H(E)\leq \dim_M(E)\leq u$.
\end{Corollary}

\begin{Corollary} Let $\{\varphi_j\}_{j=1}^{\infty}$ be infinite contractive
maps with
$$|\varphi_j(x)-\varphi_j(y)|\geq b_j |x-y|,\quad x,y\in  {\bf R}^d,\quad j=1,2,\cdots,$$
and satisfying open set condition, and denote $E$ their contractive-invariant set.
  If the equation $\sum_{j=1}^{\infty}b_j^l=1$
has finite solution $l$, then $\dim_M(E)\geq \dim_H(E)\geq \min\{d,l\}$.
\end{Corollary}

{\it Proof.} Since $\{\varphi_j\}$ are non-similar maps, we can not obtain $l^{(m)}\leq d$
from open set condition, where $l^{(m)}$ satisfies $\sum_{j=1}^mb_j^{l^{(m)}}=1$.
then similar to proof of Theorem 1, this conclusion holds.\ \ \ \ \#

\section{ Mean quadratic variations of infinite self-similar measures.}
\ \ \ \    We define
$$H(r)=\frac{1}{r^{d-\beta}}\int_{|x|\leq r}|F(x)|^2dx,$$
where $F(x)$ is the Fourier transform of $\mu$.

  If $\mu$ is a Borel measure on ${\bf R}^d$, for every $\mu$-measurable function $f$, we
use $\mu_f$ to denote the measure $\mu_f(E)=\int_Efd\mu$ for any Borel subset $E$.

 \begin{Definition}. If in addition to the definition of open set condition,
  the sets $S_j(\overline{O})$ are mutually
disjoint and $O$ intersects $K$, we call $\{S_j\}_{j\in{\bf N}}$ satisfy {\it
strong open set condition}.
\end{Definition}

  We assume $\{S_j\}_{j=1}^{\infty}$ satisfy strong open set condition. Let $d_{jk}$
denote the distance between $S_j(O)$ and $S_k(O)$ which is positive for $j\neq k$ by
strong open set condition. We assume 
\be \sum_{j\neq k}P_jP_kd_{jk}^{-\beta}<\infty.\ee
Denote $q(\lambda)=\sum_{\rho_j\leq \lambda}P_j^2\rho_j^{-\beta}$, we assume
\be q(\varepsilon\lambda)\leq\delta q(\lambda) \ee
for some $0<\varepsilon<1$ and $0<\delta<1$.

  Under the conditions (6) and (7), R.S.Strichartz [5] (P357-P358) obtained the
asymptotic property (as $r\longrightarrow +\infty$) of the function $H(r)$ and
conclusion (*) for infinite self-similar measures.

  We use $J=(j_1,j_2,\cdots, j_k)$ to denote  the multi-index, $|J|=k$ its length,
and $\Lambda$ the set of all such multi-indice, where $j_i\in {\bf N},\quad i=1,\cdots,k$
and $k\in {\bf N}$. We set 
$$P_J=P_{j_1}P_{j_2}\cdots P_{j_k}, \quad \quad \rho_J=\rho_{j_1}\cdots\rho_{j_k},\qquad E_J=E_{j_1j_2\cdots j_k}$$
For any $0<t<1$, we denote 
$$\Lambda(t)=\{J\in \Lambda:\quad \rho_J=\sup\rho_{J'},\quad \rho_{J'}<t\},$$
and for fixed parameter $\varepsilon$ (given in condition (7)), we denote
$$\Lambda_1(t)=\{J\in \Lambda (t):\quad \rho_J\geq\varepsilon t\}.$$ Then we have
\begin{Theorem} Let $\mu$ be infinite self-similar measure, we assume that the condition
(7) holds, then $V_{\beta}(t;\mu)$ is bounded below by a positive constant on $0<t\leq1$.
\end{Theorem}
{\it Proof.} Since 
$\sum_{j=1}^{\infty}P_j^2\rho_j^{-\beta}=1$, then $\sum_{J\in\Lambda(t)}P_J^2\rho_J^{-\beta}=1$.
When $J\in \Lambda_1(t)$, we have $\varepsilon t\leq\rho_J<t$. Hence
$$t^{-\beta}<\rho_J^{-\beta}\leq (\varepsilon t)^{-\beta}.$$
From the condition (7) and similar to ref.[5](P358), we can prove
$$\sum_{J\in\Lambda_1(t)}P_J^2\rho_J^{-\beta}\geq (\delta^{-1}-1)\sum_{J\in\Lambda(t)}P_J^2\rho_J^{-\beta}.$$
Hence
\begin{eqnarray*}
(\delta^{-1}-1)&=&(\delta^{-1}-1)\sum_{J\in\Lambda(t)}P_J^2\rho_J^{-\beta}\leq \sum_{J\in\Lambda_1(t)}P_J^2\rho_J^{-\beta}\\
               &\leq& \sum_{J\in\Lambda_1(t)}P_J^2(\varepsilon t)^{-\beta}\leq \sum_{J\in\Lambda(t)}P_J^2(\varepsilon t)^{-\beta},
               \end{eqnarray*}
hence $$\frac{1}{t^{\beta}} \sum_{J\in \Lambda(t)}P_J^2\geq (\delta^{-1}-1)\varepsilon^{\beta}.$$
Without loss of generality we assume $|E_0|=1$. We denote $\omega_d$ the Lebesgue
measure. Note that $\mu$ is supported by $\cup\{E_J:\quad J\in \Lambda(t)\}$ and $\mu(E_J)=P_J$.
Hence 
\begin{eqnarray*}
V_{\beta}(t;\mu)&=&\frac{1}{t^{d+\beta}}\int[\int\int\chi_{B_t(x)}(\xi)\chi_{B_t(x)}(\eta)d\mu(\xi)d\mu(\eta)]dx\\
                &=&\frac{1}{t^{d+\beta}}\int\int\omega_d(B_t(\xi)\cap B_t(\eta))d\mu(\xi)d\mu(\eta)\\
                &\geq& \frac{1}{t^{d+\beta}}\sum_{J\in\Lambda(t)}\int\int_{\xi,\eta\in E_J}\omega_d(B_t(\xi)\cap B_t(\eta))d\mu(\xi)d\mu(\eta).
                \end{eqnarray*}
Since $|E_J|=\rho_J\leq t$, hence $B_t(\xi)\cap B_t(\eta)$ contains a ball of radius
$t/2$ whenever $\xi,\eta\in E_J$. It follows that 
\begin{eqnarray*}
V_{\beta}(t;\mu)&\geq& \frac{c}{t^{\beta}}\sum_{J\in\Lambda(t)}\int\int_{\xi,\eta\in E_J}d\mu(\xi)d\mu(\eta)\\
                &\geq& \frac{c}{t^{\beta}}\sum_{J\in\Lambda(t)}P_J^2\geq c(\delta^{-1}-1)\varepsilon^{\beta},
                \end{eqnarray*}
where $c(\delta^{-1}-1)\varepsilon^{\beta}$ is a positive constant. \ \ \ \#

  From the asymptotic property of $H(r)$ of infinite self-similar measure ([5]), Theorem 4.10 and Corollary 4.12 of [2] and our
  Theorem 2, we have
\begin{Theorem} Let $\mu$ be infinite self-similar measure. Assume conditions
(6) and (7) hold, then $$\lim_{t\longrightarrow 0}(V_{\beta}(t;\mu)-P(t))=0$$
for some $P>0$ such that the following holds.

  (i) If $\{-\ln\rho_j:\quad j\in {\bf N}\}$ is non-arithematic, then $P(t)=c'$ for some
constant $c'$.

  (ii) Otherwise, let $((\ln\rho){\bf Z}),\quad\rho>1$ be the lattice generated by $\{-\ln\rho_j:\quad j\in {\bf N}\}$,
then $P(\rho t)=P(t)$.
\end{Theorem}

  From the conclusion (*) of infinite self-similar measure ([5]), Theorem 4.10 and
Corollary 4.12 of [2], if the equation $\sum_{j=1}^{\infty}\rho_j^s=1$ has finite solution $s$, then
\begin{Theorem} Let $\mu$ be infinite self-similar measure with natural weights $P_j=\rho_j^{\beta}$, where $\beta=s$
 is the finite solution of equation $\sum_{j=1}^{\infty}\rho_j^{s}=1$,
 we assume conditions
(6) and (7) holds, then for any $f\in L^2(d\mu)$ we have
$$\lim_{t\longrightarrow 0}[\frac{1}{t^{d+\beta}}\int|\mu_f(B_t(x))|^2dx-P(t)\int|f|^2d\mu]=0,$$
where $P$ defined in Theorem 3.
\end{Theorem}

\end{document}